\documentclass[12pt,leqno]{article}
\usepackage{amsmath, amsthm, amsfonts, amssymb, color}
\usepackage{mathrsfs}
\theoremstyle{definition}

\newcommand{\scr}[1]{\mathscr #1}
\definecolor{wco}{rgb}{0.5,0.2,0.3}

\numberwithin{equation}{section} \theoremstyle{remark}

\newcommand{\ua}{\uparrow}

\title{{\bf Gradient Estimate for Ornstein-Uhlenbeck Jump Processes}\footnote{Supported in
 part by WIMCS, NNSFC(10721091).}
}
\author{
{\bf Feng-Yu Wang}\\
\footnotesize{School of Mathematical Sci. and Lab. Math. Com. Sys.,
Beijing Normal
University, Beijing 100875, China}\\
\footnotesize{and}\\ \footnotesize{Department of Mathematics,
Swansea University, Singleton Park, SA2 8PP, UK}\\
\footnotesize{Email: wangfy@bnu.edu.cn; F.Y.Wang@swansea.ac.uk}}
\begin{document}
\def\R{\mathbb R}  \def\ff{\frac} \def\ss{\sqrt} \def\B{\mathbf
B}
\def\N{\mathbb N} \def\kk{\kappa} \def\m{{\bf m}}
\def\dd{\delta} \def\DD{\Delta} \def\vv{\varepsilon} \def\rr{\rho}
\def\<{\langle} \def\>{\rangle} \def\GG{\Gamma} \def\gg{\gamma}
  \def\nn{\nabla} \def\pp{\partial} \def\EE{\scr E}
\def\d{\text{\rm{d}}} \def\bb{\beta} \def\aa{\alpha} \def\D{\scr D}
  \def\si{\sigma} \def\ess{\text{\rm{ess}}}
\def\beg{\begin} \def\beq{\begin{equation}}  \def\F{\scr F}
\def\Ric{\text{\rm{Ric}}} \def\Hess{\text{\rm{Hess}}}
\def\e{\text{\rm{e}}} \def\ua{\underline a} \def\OO{\Omega}  \def\oo{\omega}
 \def\tt{\tilde} \def\Ric{\text{\rm{Ric}}}
\def\cut{\text{\rm{cut}}} \def\P{\mathbb P} \def\ifn{I_n(f^{\bigotimes n})}
\def\C{\scr C}      \def\aaa{\mathbf{r}}     \def\r{r}
\def\gap{\text{\rm{gap}}} \def\prr{\pi_{{\bf m},\varrho}}  \def\r{\mathbf r}
\def\Z{\mathbb Z} \def\vrr{\varrho} \def\ll{\lambda}
\def\L{\scr L}\def\Tt{\tt} \def\TT{\tt}\def\II{\mathbb I}
\def\i{{\rm in}}\def\Sect{{\rm Sect}}\def\E{\mathbb E} \def\H{\mathbb H}
\def\M{\scr M}\def\Q{\mathbb Q} \def\texto{\text{o}} \def\LL{\Lambda}
\def\Rank{{\rm Rank}} \def\B{\mathbf B} \def\i{{\rm i}}

\maketitle
\begin{abstract} By using absolutely continuous lower bounds of the L\'evy measure, explicit gradient estimates are derived for the semigroup of the corresponding L\'evy process with a linear drift. A derivative formula  is presented for the conditional distribution of the process at time $t$ under the condition that the process jumps before $t$. Finally, by using bounded perturbations of the L\'evy measure, the resulting gradient estimates are extended to linear SDEs driven by L\'evy-type processes.
\end{abstract} \noindent

 \vskip 2cm

\section{Introduction}

It is well-known that a L\'evy process can be decomposed into two
independent parts, i.e. the diffusion  part and the jump part. If the
diffusion part is non-degenerate, regularity properties for the
semigroup of the Brownian motion can be easily confirmed for the
L\'evy semigroup.  On the other hand, when the L\'evy process is pure jump, existence and regularities of the transition density have been
derived by using conditions on the symbol or the L\'evy measure (see \cite{KS,KS2,PZ} and references within);
see also \cite{BJ, JS} for heat kernel upper bounds for $\aa$-stable processes
with drifts. As a continuation to the recent work
\cite{W10}, where the  coupling property
and applications are studied by using absolutely continuous lower bounds of the  L\'evy measure,   this note aims to derive gradient estimates
of the L\'evy semigroup in the same spirit.

Let $L_t$ be the L\'evy process on $\R^d$ with symbol (see e.g.
\cite{A})

$$\eta(u)= \i \<u,b\> -\<Q u,u\> +\int_{\R^d} \big(\e^{\i\<u,z\>}-1-\i\<u,z\>1_{\{|z|<1\}}\big)\nu(\d z),$$
where $b\in \R^d$, $Q$ is a non-negatively definite $d\times d$
matrix, and $\nu$ is a L\'evy measure on $\R^d$. In references the L\'evy symbol is also called the characteristic exponent or the L\'evy exponent, and in e.g. \cite{J}, $-\eta$ rather than $\eta$ is called the L\'evy symbol.
It is well known that  $L_t$ is a strong Markov process on $\R^d$ generated by

\beq\label{G1} \scr L f := \<b,\nn f\>+ \text{Tr}(Q\nn^2 f)+\int_{\R^d} \big\{f(z+\cdot)-f
-\<\nn f, z\>1_{\{|z|\le 1\}}\big\}\nu(\d z)\end{equation} for $ f\in C_b^2(\R^d).$

Let $P_t$ be the semigroup for the solution of the linear stochastic
differential equation

\beq\label{E1}\d X_t= AX_t\d t +\d L_t,\end{equation} where $A$ is a $d\times d$ matrix. According to \cite{BRS}, we have

\beq\label{M} P_t f(x)= \int_{\R^d} f(\e^{tA}x+y)\mu_t(\d
y),\end{equation} where $\mu_t$ is the probability measure on $\R^d$
with characteristic function

\beq\label{M2} \hat \mu_t(z)= \exp\bigg[\int_0^t \eta(\e^{s A^*}
z)\d s\bigg],\ \ \ z\in\R^d.\end{equation}   Let $\scr B_b(\R^d)$ be the set of all
bounded measurable functions on $\R^d$.
We shall estimate $\|\nn P_t f\|_\infty$,  the uniform norm of the
gradient $\nn P_tf$, for $t>0$ and   $ f\in\scr B_b(\R^d).$
When the L\'evy measure is finite, with a positive probability
the process does not jump before a fixed time $t>0$. So, in this case, the semigroup is not strong Feller and thus, does not have finite uniform gradient estimate. Therefore,  to
derive the uniform gradient estimate, it is
essential to assume that $\nu$ is infinite. Since $\nu$ is always finite
outside a neighborhood of $0$, the behavior of $\nu$ around the
origin will be crucial for the study.

We will make use of the following lower bound condition of $\nu$:

\beq\label{C1} \nu(\d z)\ge |z|^{-d}S(|z|^{-2})1_{\{|z|<r_0\}}\d
z,\end{equation}  where $r_0\in (0,\infty]$ is a constant and $S$ is a Bernstein function with $S(0)=0.$ Let

\beg{equation*}\beg{split} &c_0=  \int_{\{|z|\le \e^{-\|A\|}\}} (1-\cos z_1) |z|^{-d}\d z,\\
&\ll_0= \int_{\R^d} (r_0\lor |z|)^{-d} S((r_0\lor |z|)^{-2}) \d z,\end{split}\end{equation*} where $z_1$ stands for the first coordinate of $z$, and $\|A\|$ is the operator norm of $A$. We have
$c_0\in (0,\infty).$ Since $S(r)\le cr$ holds for some constant $c\in (0,\infty)$, we have  $\ll_0<\infty.$  In particular, if $r_0=\infty$ then $\ll_0=0.$ We will estimate $\|\nn P_t f\|_\infty$ by using the upper bound of $A$ and the function
$$\aa(t):= \int_0^\infty \ff 1 {\ss r}\e^{-t S(r)}\d r,\ \ \ t>0.$$
Obviously, if $\lim_{r\to\infty} \ff{S(r)}{\log r} =\infty$ then $\aa(t)<\infty$ for all $t>0.$

\beg{thm}\label{T1.1} Let $(\ref{C1})$ hold and let $c_0,\ll_0, \aa(t)$ be defined above, let $\theta\in \R$ be such that $A\le -\theta I$. Then there exists a constant $c_1\in (0,\infty)$ depending only on $d$ and $\theta$ such that

\beq\label{G}\|\nn P_t f\|_\infty\le \|f\|_\infty c_1 \e^{\ll_0(t\land 1)-\theta^+t}\bigg\{ \aa(c_0(t\land 1))
+\ff {(t\land 1)S(r_0^{-2}) } {r_0 }\bigg\} \end{equation} holds for any $t>0$ and $f\in \scr B_b(\R^d).$
If moreover $A=0$, then there exists $c_1$ depending on $d$ such that

\beq\label{G2} \|\nn P_t f\|_\infty \le \|f\|_\infty \e^{\ll_0t}\bigg\{ \ff 1 {\ss {2\pi}} \aa(c_0 t) + \ff{c_1 (1-\e^{-t\ll_0})S(r_0^{-2})}{r_0 \ll_0}\bigg\} \end{equation} holds for any $t>0$ and $f\in \scr B_b(\R^d),$ where
$\ll_0=\ff{1-\e^{-t\ll_0}}{r_0\ll_0}=0$ for $r_0=\infty$. \end{thm}

Now, we consider the gradient estimate for the semigroup associated to the linear SDE driven by a L\'evy-type process. Let $\si(x,\d y)$ be a
signed kernel on $\R^d$, i.e. for each $x\in \R^d$, $\si(x,\cdot)$ is a signed measure while for each measurable set $A$, $\si(\cdot, A)$
is a measurable function. We call $\si$ bounded if

$$\|\si\|_\infty:= \sup_{x\in \R^d} |\si(x,\cdot)|(\R^d)<\infty.$$ Let $L_t^{+\si}$ be the L\'evy-type process with jump measure

$$q(x,\d z):= \nu(\d z-x)+ \si(x,\d z)$$ for a bounded $\si$. In other words, there exist $b\in \R^d$ and non-negatively definite $d\times d$-matrix $Q$ such that $L_t^{+\si}$ is generated by

\beq\label{GG} \scr L^{+\si} f(x) = \scr Lf(x)+ \int_{\R^d} \big\{f(z)-f(x)\big\}\si(x,\d z)=:\scr L f(x)+\si f(x)
\end{equation} for $f\in C^2_b(\R^d)$, where $\scr L$ is in (\ref{G}).
Let $P_t^{+\si}$ be the semigroup associated to the linear SDE

$$\d X_t = A X_t\d t+ \d L_t^{+\si}.$$ Combining Theorem \ref{T1.1} with a standard perturbation argument, we prove the following result on the gradient estimate of $P_t^{+\si}$.

\beg{cor}\label{C1.3}  If $(\ref{C1})$ holds for some $S$ such that $\int_0^1 \aa(t)\d t<\infty$, then there
 exists a constant $c\in (0,\infty)$ such that

 $$\|\nn P_t^{+\si}f\|_\infty\le c \big\{\aa(c_0(t\land 1)) + \|\si\|_\infty\big\}\|f\|_\infty,\ \ t>0, f\in\scr B_b(\R^d)$$ holds for any bounded $\si$.
 \end{cor}

To illustrate our results, we consider below two typical choices of
$\rr_0$.

\beg{exa} \label{E1.3} $(1)$ If $\nu(\d z)\ge c |z|^{-d-\aa}1_{\{|z|\le r_0\}}$
for some $c,r_0>0$ and $\aa\in (0,2)$, then

$$\|\nn P_tf\|_\infty \le \ff{c'}{(t\land 1)^{1/\aa}}
\e^{-\theta^+t}\| f\|_\infty,\ \ t>0, f\in \scr B_b(\R^d)$$ holds for some constant $c'\in (0,\infty).$
If $\aa\in (1,2)$, then there exists a constant $c\in (0,\infty)$ such that

$$\|\nn P_t^{+\si}f\|_\infty \le  c \|f\|_\infty\Big\{\ff 1 {(t\land 1)^{1/\aa}} +\|\si\|_\infty\Big\},\ \
t>0, f\in \scr B_b(\R^d)$$ holds for any bounded $\si$.

$(2)$ If   $\nu(\d z)\ge c
|z|^{-d}\log^{1+\vv}(1+|z|^{-2})1_{\{|z|\le r_0\}}$ for some
$c,r_0,\vv>0$, then

$$\|\nn P_tf\|_\infty \le c_1\| f\|_\infty
\exp[c_2t^{-1/\vv} -\theta^+t], \ \ t>0, f\in\scr B_b(\R^d)$$ holds for some constants
$c_1,c_2\in(0,\infty).$\end{exa}

 Since  for the $\aa$-stable process one has (see  Corollary \ref{C2.2} (2) below for a more general result)
 $$\sup_{\|f\|_\infty\le 1}\|P_tf\|_\infty \ge \ff c {t^{1/\aa}}$$ for some constant $c>0$. Thus, the upper bound in
 Example \ref{E1.3}(1) is sharp.

\

The main idea of the proof is to compare the process with the
$S$-subordinate semigroup of the Brownian motion. To this end, we
shall study in the next section the gradient estimate for
subordinate semigroups. We will see that to compare the original
semigroup with the subordinate semigroup, the error term is given by the conditional
distribution of a compound Poisson process under the condition that the process jumps before time $t$.   Thus, in Section 3 we will
study the gradient estimate for the corresponding conditional distribution for  compound Poisson processes. In this case, a   derivative formula is presented. By combining
results derived in Section 2 and Section 3, we prove Theorem \ref{T1.1} in Section 4.  Finally, the proofs of
Corollaries \ref{C1.3}  and Example \ref{E1.3} are addressed in Section 5.

\section{Gradient estimates for subordinate semigroups}

This section is a counterpart of the recent work \cite{GRW} where
dimension-free Harnack inequality is investigated for subordinate
semigroups, see e.g. \cite{SV} and references within for potential theory and historical
remarks on subordinations of the Brownian motion.

Let $(E,\rr)$ be a Polish space. For a function $f$ on $E$, define

$$|\nn f|(x):= \limsup_{y\to x}\ff{|f(y)-f(x)|}{\rr(x,y)},\ \ x\in
E.$$ Let $P_t^0$ be a (sub-)Markov semigroup on $\scr B_b(E)$ such
that for some positive function $\varphi$ on $(0,\infty)$,

\beq\label{2.1} |\nn P_t^0 f|\le \|f\|_\infty \varphi(t),\ \ t>0, f\in
\scr B_b(E)\end{equation} holds. We intend to estimate the gradient of a
subordinate semigroup $P_t^S$ of $P_t^0$ induced by a Bernstein
function $S$. More precisely, for any $t\ge 0$ let $\mu_t^S$ be the
probability measure on $[0,\infty)$ with Laplace transformation

\beq\label{2.3} \int_0^\infty \e^{-\ll s}\mu_t^S(\d s)= \e^{-t
S(\ll)},\ \ \ll\ge 0.\end{equation} Then the $S$-subordination of $P_t^0$ is  given  by

\beq\label{2.2} P_t^S=\int_0^\infty P_s^0\mu_t^S(\d s),\ \ t\ge
0.\end{equation}
The following assertion follows immediately from (\ref{2.2}) and the
dominated convergence theorem.

\beg{thm}\label{T2.1} If $(\ref{2.1})$ holds with $\int_0^\infty
\varphi(s) \mu_t^S(\d s)<\infty$, then

$$|\nn P_t^Sf|\le \|f\|_\infty\int_0^\infty \varphi(s)\mu_t^S(\d
s),\ \ f\in \scr B_b(E).$$ \end{thm}

In particular, we have  the following explicit gradient estimates  by using known results on diffusion semigroups.

\beg{cor}\label{C2.2} $(1)$ Let $E$ be a complete connected Riemannian
manifold and $P_t^0$ be the diffusion semigroup generated by $\DD+Z$
for a vector field $Z$ on $E$ such that

$$\Ric -\nn Z\ge 0$$ holds. Then

$$\|\nn P_t^S f\|_\infty\le \ff{ \|f\|_\infty}{\ss{2\pi}} \int_0^\infty \ff{1}{\ss r}\, \e^{-tS(r)}\d r,\ \
t>0, f\in \scr B_b(E).$$

$(2)$ Let $P_t^0$ be generated by $\DD$ on $\R^d$. We have

$$\sup_{\|f\|_\infty\le 1} \|\nn P_t^Sf\|_\infty\ge \ff 1 {\ss 2\,\pi} \int_0^\infty \ff 1 {\ss r} \e^{-tS(r)}\d r.$$ \end{cor}

\beg{proof} (1) It is well-known that the curvature condition implies
(cf. \cite{B})

$$P_t^0 f^2-(P_t^0f)^2 \ge  t  |\nn P_t^0 f|^2.$$ This implies that

$$\|\nn P_t^0f\|_\infty \le \ff 1 {\ss t} \|f\|_\infty.$$ Then the
proof of (1) is  finished by combining this with Theorem
\ref{T2.1} and noting that

\beg{equation*}\beg{split} &\int_{0}^{\infty}\frac{\mu_{t}^{S}\left(
\d s\right)}{\ss s}=\int_{0}^{\infty}\frac{1}{\ss{2\pi}}\int_{0}^{\infty} \ff 1 {\ss r} \e^{-rs}\d r \mu_{t}^S\left( \d s\right)
\\
&=\frac{1}{\ss{2\pi}}\int_{0}^{\infty}\ff 1 {\ss r} \int_{0}^{\infty}
\e^{-rs}\mu_{t}^S\left( \d s\right)\d r=\frac{1}{\ss{2\pi}}\int_{0}^{\infty}r^{-1/2} \e^{-tS\left(r\right)}\d r.\end{split}\end{equation*}

(2) Let $P_t^0$ be generated by $\DD$ on $\R^d$. We have

$$P_s^0 f(x) = \ff 1 {(4\pi s)^{d/2}} \int_{\R^d} \e^{-|x-y|^2/(4s)}f(y)\d y.$$ Take

$$f(x)= 1_{[0,\infty)}(x_1)- 1_{(-\infty,0)}(x_1).$$ We have $\|f\|_\infty =1$ and

\beg{equation*}\beg{split} P_s^0 f(x) &=\ff 1 {2\ss{\pi s}} \bigg\{\int_0^\infty \e^{-(r-x_1)^2/(4s)}\d r -
\int_{-\infty}^0  \e^{-(r-x_1)^2/(4s)}\d r\bigg\}\\
&= \ff 1 {2\ss{\pi s}} \bigg\{\int_{-x_1}^\infty \e^{-r^2/(4s)}\d r -
\int_{-\infty}^{-x_1}  \e^{-r^2/(4s)}\d r\bigg\}.\end{split}\end{equation*}
So,

$$\ff{\d}{\d x_1}P_s^0 f(x)= \ff 1 {\ss{\pi s}} \e^{-x_1^2/(4s)}\le \ff 1 {\ss{\pi s}},\ \ s>0, x\in\R^d.$$ Combining this with (\ref{2.2}) and using the dominated convergence theorem, we arrive at

$$\ff{\d}{\d x_1} P_t^Sf(x)\Big|_{x=0} = \ff 1 {\ss\pi}\int_0^\infty\ff 1 {\ss s} \mu_t^S(\d s)= \ff 1 {\ss 2\, \pi}
\int_0^\infty \ff 1 {\ss r }\e^{-tS(r)}\d r.$$
\end{proof}

\section{A derivative formula}

Let $\nu(\d z)\ge \rr_0(z)\d z=:\nu_0(\d z)$ for some non-negative measurable
function $\rr_0$ on $\R^d$ such that

\beq\label{3.1} \ll_0:= \int_{\R^d} \rr_0(z) \d z\in
(0,\infty).\end{equation}  Let $(L_t^0)_{t\ge 0}$ be the compound
Poisson process with L\'evy measure $\nu_0$. Then $L_t^0$ can be
realized as

\beq\label{3.2} L_t^0= \sum_{i=1}^{N_t} \xi_i,\ \ \ t\ge
0,\end{equation} where $N_t$ is the Poisson process with rate
$\ll_0$ and $\{\xi_i\}$ are i.i.d. random variables on $\R^d$ which are independent of
$(N_t)_{t\ge 0}$ and have  common distribution $\nu_0/\ll_0.$ Here, we
set $\sum_{i=1}^0\xi_i=0$ by convention. Let $(L_t^1)_{t\ge 0}$  be the
L\'evy process which is  independent of $(L_t^0)_{t\ge 0}$ and has
L\'evy measure $\nu-\nu_0$, such that

\beq\label{3.3} L_t:= L_t^1+L_t^0,\ \ t\ge 0\end{equation} is the L\'evy process with
symbol $\eta$.  As we explained in the Introduction,
to ensure the strong Feller property for a jump process, it is
essential to restrict on the event that the process jumps before a
fixed time. Thus, instead of $P_t$, it is natural for us to  investigate the gradient estimate for
$P_t^1$ defined by
$$ P_t^1 f(x)= \E\Big\{f(X_t^x)1_{\{N_t\ge 1\}}\Big\},\ \ f\in \scr
B_b(\R^d),\ t>0,$$ where $X_t^x$ solves (\ref{E1}) with initial data $x$. The following result provides a derivative
formula for this operator, which can be regarded as the jump counterpart of the
Bismut-Elworthy-Li   formula for diffusion processes
\cite{Bismut, EL}.

\beg{thm}\label{T3.1} Let $\rr_0$ be non-negative and differentiable such that $\nu(\d z)\ge \rr_0(z)\d z$,
$\ll_0:=\int_{\R^d} \rr_0(z)\d z\in (0,\infty)$, and

\beq\label{3.5} \int_{\R^d} \Big\{\sup_{x:|x-z|\le \vv} |\nn \rr_0|(x)\Big\}\d
z<\infty\end{equation} holds for some $\vv>0$. Then for any $t>0$
and $f\in \scr B_b(\R^d),$

\beq\label{3.6} \nn P_t^1 f(x)= -\E\Big\{f(X_t^x)1_{\{N_t\ge 1\}}\ff
1 {N_t} \sum_{i=1}^{N_t} \e^{A^*\tau_i}\nn\log
\rr_0(\xi_i)\Big\},\end{equation} where $\tau_i$ is the $i$-th jump
time of $(N_t)_{t\ge 0}$ and $A^*$ is the transposition of $A$.
Consequently, if $A\le -\theta I$ then

$$\|\nn P_t^1 f\|_\infty \le \|f\|_\infty\ff{\e^{\theta^- t}(1-\e^{-\ll_0 t})}{\ll_0}
\int_{\R^d} |\nn\rr_0|(z)\d z,\ \ t>0, f\in \scr
B_b(\R^d).$$\end{thm}

\beg{proof} We shall make use of a formula for random shifts of the
compound Poisson process derived in \cite{W10}. Let $\LL(\d w)$ be the
distribution of $L^0:=(L_t^0)_{t\ge 0}$ which is a probability measure on the path space

 $$W= \Big\{\sum_{i=1}^\infty x_i1_{[t_i,\infty)}: i\in \N, x_i\in \R^d\setminus\{0\}, \ 0\le t_i\uparrow \infty\ \text{as}
 \ i\uparrow \infty\Big\}$$ equipped with the $\si$-algebra induced by $\{ w\mapsto  w_t:\ t\ge 0\}.$

  Let  $(\tau,\xi)$ be
a $[0,t]\times \R^d$-valued random variable such that the joint
distribution of $(L^0, \tau,\xi)$ is

$$g( w, s,z)\LL(\d w)\d s\nu_0(\d z).$$ Let
$\DD w_t= w_t- w_{t-}$ and

$$U( w)= \sum_{\DD w_t\ne0} g( w-\DD w_t 1_{[t,\infty)}, t, \DD w_t).$$
By \cite[Corollary 2.3]{W10}, for any bounded measurable function
$F$ on the path space of $L^0$, one has

\beq\label{T} \E\big(F1_{\{U>0\}}\big)(L^0) =
\E\Big\{\ff{F1_{\{U>0\}}}{U}\Big\}(L^0+\xi1_{[\tau,\infty)}).\end{equation}
Now, let $(\tau,\xi)$ be independent  of $(L_t^1, L_t^0)_{t\ge 0}$
with distribution

$$\ff 1 {t\ll_0} 1_{[0,t]}(s) \d s \nu_0(\d z).$$  We have  $g( w, s,z)= \ff1 {\ll_0 t} 1_{[0,t]}(s).$ Since
$\tau$ is independent of $L^0$ so that with probability one $\tau (\le t)$ is not a jump time of $L^0$, and since
$\xi\ne 0$ a.s., we have

$$ U(L^0+\xi1_{[\tau,\infty)})= \ff{N_t+1}{\ll_0 t}.$$
Since $Y_t:= \int_0^t\e^{(t-s)A} \d L_s^1$   is independent of

$$  \e^{At}x + \int_0^t\e^{A(t-s)}\d L_s^0,$$ it follows from
(\ref{T}) that for any $z_0\in\R^d$ and $\vv\in (-1,1)$,

\beq\label{3.8} \beg{split}&P^1_t f(x+\vv z_0)=    \E\bigg\{
f\bigg(Y_t + \e^{At}(x+\vv z_0)+\int_0^t\e^{A(t-s)}\d
L_s^0\bigg)1_{\{N_t\ge
1\}}\bigg\}\\
&=\ll_0t\E\bigg\{\ff{ f\big(Y_t + \e^{At}(x+\vv
z_0)+\int_0^t\e^{A(t-s)}\d\{
L^0+\xi1_{[\tau,\infty)}\}_s\big)}{N_t+1}\bigg\}\\
&=\ll_0t\E\bigg\{ \ff{f\big(Y_t + \e^{At}x   +\int_0^t\e^{A(t-s)}\d\{
L^0+(\xi+\vv\e^{A\tau}
z_0)1_{[\tau,\infty)}\}_s\big)}{N_t+1}\bigg\}.\end{split}\end{equation}
On the other hand, since the joint distribution of $(L^0, \tau,
\xi+\vv \e^{A\tau}z_0)$ is

$$\ff 1 {\ll_0 t} 1_{[0,t]}(s)
\ff{\rr_0(z-\vv\e^{As}z_0)}{\rr_0(z)}\LL(\d w)\d s\nu_0(\d z),$$
(\ref{T}) holds for $\xi':=\xi+\vv\e^{\tau A}z_0$ in place of $\xi$
with

$$U(L^0)= \ff 1 {\ll_0 t} \sum_{i=1}^{N_t}
\ff{\rr_0(\xi_i-\vv \e^{\tau_i A}z_0)}{\rr_0(\xi_i)}.$$ Consequently, for any $F\ge 0$, using $FU$ in place of $F$ in (\ref{T}) one obtains

$$\E\big\{F(L^0)U(L^0) 1_{\{N_t\ge 1\}}\big\}= \E F(L^0 +\xi' 1_{[\tau,\infty)}).$$ Taking $n_t(w)= \sum_{s\le t} 1_{\{\DD w_s\ne 0\}}$ and

$$F( w)= \ff{f\big(z+ \int_0^t\e^{(t-s)A} \d w_s\big)}{n_t( w)}1_{\{n_t(w)\ge 1\}},\ \ w\in W$$ for $z\in\R^d$, we arrive at

\beg{equation*}\beg{split} &\ff 1 {\ll_0 t}\E\bigg\{ f\bigg(z+ \int_0^t\e^{(t-s)A} \d L_s^0\bigg)
 \ff{1_{\{N_t\ge 1\}}}{N_t} \sum_{i=1}^{N_t} \ff{\rr_0(\xi_i-\vv \e^{A\tau_i} z_0)}{\rr_0(\xi_i)}\bigg\}\\
 &= \E\bigg\{ \ff{f\big(z   +\int_0^t\e^{A(t-s)}\d\{
L^0+(\xi+\vv\e^{A\tau}
z_0)1_{[\tau,\infty)}\}_s\big)}{N_t+1}\bigg\},\ \ z\in\R^d.\end{split}\end{equation*}
Combining this with  (\ref{3.8}), we obtain

$$P^1_t f(x+\vv z_0)= \E\bigg\{f(X_t^x)1_{\{N_t\ge 1\}} \ff 1
{N_t}\sum_{i=1}^{N_t}
\ff{\rr_0(\xi_i-\vv\e^{A\tau_i}z_0)}{\rr_0(\xi_i)}\bigg\}.$$ Therefore,  for any $\vv\ne 0$ we have

\beq\label{AA}\ff{P_t^1f(x+\vv z_0) -P_t^1f(x)}\vv=
\E\Big\{f(X_t^x)1_{\{N_t\ge 1\}} \ff 1 {N_t}\sum_{i=1}^{N_t}
\ff{\rr_0(\xi_i-\vv
\e^{A\tau_i}z_0)-\rr_0(\xi_i)}{\vv \rr_0(\xi_i)}\Big\}.\end{equation}
Noting that for $i\le N_t$ one has $\tau_i\le t$ so that
$\e^{A\tau_i} z_0$ is bounded, and noting that for each $i$ one has
$$\lim_{\vv\downarrow 0} \ff{\rr_0(\xi_i-\vv
\e^{A\tau_i}z_0)-\rr_0(\xi_i)}{\vv \rr_0(\xi_i)}= -\<\e^{A\tau_i} z_0,\nn\log \rr_0(\xi_i)\>
=-\<z_0, \e^{A^*\tau_i}\nn\log \rr_0(\xi_i)\>,$$
by (\ref{3.5}) we are able to use the
dominated convergence theorem to derive (\ref{3.6}) by letting
$\vv\to 0$ in (\ref{AA}). \end{proof}

\section{Proof of Theorem \ref{T1.1}}

\subsection{Proof of (\ref{G2}) for $A=0$} We shall first consider the case where
$r_0=\infty$ then pass to finite $r_0$ by using Theorem \ref{T3.1}.

\

(I) For $r_0=\infty$, i.e.
\beq\label{R} \nu(\d z)\ge |z|^{-d}S(|z|^{-2})\d z.\end{equation}
Then  \beg{equation*}\beg{split} &\eta_1(u) := \int_{\R^d}
\big(\e^{\i\<u,z\>}-1-\i\<u,z\>
1_{\{|z|<1\}}\big)|z|^{-d}S(|z|^{-2})\d z\\
&\eta_2(u):=\eta(u)-\eta_1(u)\\
&\qquad=\i \<u,b\> -\<Q u,u\> +\int_{\R^d}
\big(\e^{\i\<u,z\>}-1-\i\<u,z\>1_{\{|z|<1\}}\big)\big\{\nu(\d
z)-|z|^{-d}S(|u|^2)\d z\big\}\end{split}\end{equation*} provide two L\'evy
symbols. Noting that $S(|z|^{-2})\ge 1_{\{|z|\le |u|^{-1}\}}
S(|u|^2)$ and

\beg{equation*}\beg{split}&-\int_{\{|z|\le |u|^{-1}\}}
\big(\e^{\i\<u,z\>}-1-\i\<u,z\> 1_{\{|z|<1\}}\big)|z|^{-d}\d z\\
&= \int_{\{|z|\le |u|^{-1}\}}(1-\cos \<u,z\>)|z|^{-d}\d
z\\
&= \int_{\{|z|\le 1\}}\Big(1-\cos \Big\<\ff u
{|u|},z\Big\>\Big)|z|^{-d}\d z\\
&= \int_{\{|z|\le 1\}}(1-\cos z_1)|z|^{-d}\d
z=c_0\in (0,\infty),\end{split}\end{equation*} we see that

\beg{equation*}\beg{split}&u\mapsto \eta(u)+ c_0
S(|u|^2)\\
&=\eta_2(u)+ \int_{\R^d} \big(\e^{\i\<u,z\>}-1-\i\<u,z\>
1_{\{|z|<1\}}\big)|z|^{-d}\big\{S(|z|^{-2})-S(|u|^2)1_{\{|z|\le
|u|^{-1}\}}\big\}\d z\end{split}\end{equation*} is also a L\'evy symbol.
Let $P_t^S$ be the semigroup of the L\'evy process with L\'evy symbol $-c_0
S(|\cdot|^2)$, and let $\tt P_t^S$ be the one with L\'evy symbol $\eta+c_0
S(|\cdot|^2)$. We have

\beq\label{A1} P_t = P_t^S \tt P_t^S.\end{equation} Since $P_t^S$ is
the $c_0S$-subordination of the semigroup generated by $\DD$ on
$\R^d$, according to  Corollary \ref{C2.2} for $E=\R^d$ and $Z=0$,

\beq\label{A1'} \|\nn P_t^S f\|_\infty \le \|f\|_\infty
\int_0^\infty \ff 1 {\ss{2\pi r}} \e^{-c_0tS(r)}\d r =\ff 1 {\ss{
2\pi}}\aa(c_0t)\|f\|_\infty.\end{equation} Combining this with
(\ref{A1}) we derive

\beq\label{A2} \|\nn P_t f\|_\infty\le \ff 1 {\ss {2\pi}}
\aa(c_0t)\|f\|_\infty.\end{equation} Thus, the desired assertion
holds if $r_0=\infty.$

\

(II) For $r_0\in (0,\infty).$  Take

$$\rr_0(z)= (r_0\lor |z|)^{-d} S((r_0\lor |z|)^{-2}).$$Then

 \beq\label{A3}\bar\nu(\d z):= \nu(\d z)+ \rr_0(z)\d z\ge |z|^{-d}
S(|z|^{-2})\d z.\end{equation}
 Let $\bar L_t^0$ be the compound
Poisson process with L\'evy measure $\rr_0(z)\d z$, and let

$$\bar P_t^1f(x)=\E\big\{1_{\{\bar\tau_1\le t\}}f(x+\bar L_t^0)\big\},$$
where $\bar\tau_1$ is the first jump time of $\bar L_t^0$. Let $L_t$
be the L\'evy process with L\'evy symbol $\eta$ which is independent of
$\bar L_t^0$. Then $\bar L_t := L_t+\bar L_t^0$ is the L\'evy
process with L\'evy symbol

$$u\mapsto  \eta(u)+ \int_{\R^d}
(\cos\<u,z\>-1)\rr_0(z)\d z.$$ Therefore,

\beg{equation*}\beg{split} &\bar P_t f(x):= \E f(x+ \bar L_t)\\
&= \E\big\{f(x+L_t)1_{\{\bar\tau_1>t\}}\big\}+ \E\big\{f(x+L_t+\bar
L_t^0)1_{\{\bar\tau_1\le t\}}\big\}\\
&=\e^{-\ll_0 t}P_tf(x)+ \bar P_t^1P_t
f(x).\end{split}\end{equation*} This implies that

\beq\label{A4} P_t f(x) =\e^{\ll_0 t} \big(\bar P_t f- \bar P_t^1
P_t f)(x).\end{equation} According to (\ref{A3}) and (I), (\ref{A2})
holds for $\bar P_t$ in place of $P_t$, i.e.

\beq\label{A5} \|\nn \bar P_t f\|_\infty\le \ff 1 {\ss {2\pi}}
\aa(c_0t)\|f\|_\infty.\end{equation} On the other hand, we have

$$ |\nn\rr_0(z)|\le 1_{\{|z|\ge r_0\}} \big\{d |z|^{-d-1} S(r_0^{-2})
+ 2 |z|^{-d-3} S'(|z|^{-2})\big\}.$$ Since  $S'$ is decreasing, $S$ is increasing and $S(0)=0$, from this we may find a constant  
 $c $ depending only on $d$ such that

\beg{equation*}\beg{split} &\int_{\R^d}\Big\{\sup_{x: |x-z|< r_0/2}|\nn \rr_0(x)|\Big\}\d z \le c \int_{r_0}^\infty
r^{-2}\big\{S(r_0^{-2})+ r^{-2} S'(r^{-2}/4)\big\}\d r\\
&=c\int_{r_0}^\infty\Big\{\ff{S(r_0^{-2})}{r^2} -\ff 2 r \ff{\d}{\d r} S(r^{-2}/4)\Big\}\d r\le  \ff c{ r_0} S(r_0^{-2}) + \ff{2 c} {r_0} S(r_0^{-2}/4) 
 \le \ff{3c}{r_0} S(r_0^{-2}).\end{split}\end{equation*} Therefore, it follows from    Theorem \ref{T3.1} with  $\theta=0$ that
 
\beq\label{DDD}\beg{split}\|\nn \bar P_t^1f\|_\infty &\le \ff{3cS(r_0^{-2}) (1-\e^{-\ll_0 t})}{r_0 \ll_0} \|f\|_\infty\\
&\le
\ff{3cS(r_0^{-2})  t}{r_0 } \|f\|_\infty,\ \ t>0.\end{split}\end{equation}
 Combining this with (\ref{A4}) and
(\ref{A5}) we obtain the desired gradient estimate (\ref{G2}).

\subsection{Proof of (\ref{G}) for $A\ne 0$}
(III)  We first observe
that it suffices to prove (\ref{G})  for $t\in (0,1].$  
  Assume that (\ref{G}) holds for $t\in (0,1]$. By the semigroup
property we have 

$$|\nn P_t f|\le |\nn P_{t\land 1}(P_{(t-1)^+}f)|\le c_1
 \aa(c_0(t\land 1))\|f\|_\infty,\ \ t>0 $$   for some constant
$c_0,c_1\in (0,\infty).$ So, the desired inequality (\ref{G}) holds for
$\theta\le 0.$ Next,  since $A\le -\theta
I$ implies that $|X_t^x-X_t^y|\le \e^{-\theta t}
|x-y|,$      we
have

\beg{equation*} \beg{split}\ff{|P_tf(x)-P_tf(y)|}{|x-y|} &\le \ff
{|\E P_1
f(X_{t-1}^x)- \E P_1f(X_{t-1}^y)|}{|x-y|}\\
&\le \e^{-\theta (t-1)} \E\bigg\{\ff{  | P_1 f(X_{t-1}^x)-
P_1f(X_{t-1}^y)|}{|X_{t-1}^x-X_{t-1}^y|}\bigg\}.\end{split}\end{equation*}
Letting $y\to x$ and using the assertion for $t=1$ and the dominated
convergence theorem, we arrive at

$$|\nn P_t f(x)|\le \e^{-\theta (t-1)}|\nn P_1 f(X_{t-1}^x)|\le c_1
\e^{-\theta (t-1)} \aa(c_0(t\land 1))\|f\|_\infty,\ \ t>1.$$ That
is, (\ref{G}) holds also for $t>1$ with a different constant $c_1$.

\ \newline  (IV) For $r_0=\infty$ and $t\in (0,1].$
Let

\beg{equation*}\beg{split} \eta_1(u) &= \int_{\R^d} \big(\e^{\i\<u,z\>}-1-\i\<u,z\>
1_{\{|z|<1\}}\big)|z|^{-d}S(|z|^{-2})\d z\\
&= \int_{\R^d}\big(\cos \<u,z\> -1\big)|z|^{-d}S(|z|^{-2})\d
z,\end{split}\end{equation*}  and  $\eta_2 =\eta-\eta_1.$ By
(\ref{R}), both $\eta_1$ and $\eta_2$ are L\'evy symbols.  We have

\beg{equation*}\beg{split} & \eta_1(\e^{sA^*}u)+ c_0S(|u|^2)\\
  &= \int_{\R^d} \big(\cos \<z, \e^{sA^*}u\> -1\big)|z|^{-d}S(|z|^{-2})\d z+c_0S(|u|^2)\\
   &=\int_{\R^d}\Big(\cos \Big\<z, \ff{\e^{sA^*}u}{|\e^{s A^*}u|}\Big\> -1\Big)|z|^{-d}S(|z|^{-2}|\e^{sA^*}u|^{2})
   \d z+ c_0S(|u|^2) \\
   &= \int_{\R^d}\big(\cos z_1 -1\big)|z|^{-d}\big\{S(|z|^{-2}|\e^{sA^*}u|^{2})
   - S(|u|^2)1_{\{|z|\le \e^{-\|A\|}\}}\big\} \d z\\
   &= \int_{\R^d} \big( \e^{\i\<u,z\>}-1-\i\<u,z\>1_{\{|z|<1\}}\big) |z|^{-d}
   \big\{S(|z|^{-2}|\e^{sA^*}u|^{2})- S(|u|^2)1_{\{|z|\le \e^{-\|A\|}\}}\big\}
   \d z.\end{split}\end{equation*}  Since for $s\in
[0,1]$

$$S(|z|^{-2}|\e^{s A^*}u|^2)\ge S(|u|^2)1_{\{|u|\le \e^{-\|A\|}\}},$$ this implies  that

   $$u\mapsto \eta_1(\e^{s A^*}u) + c_0 S(|u|^2)$$ is a L\'evy symbol. In particular,
   there exists a probability measure $\pi_t$ on $\R^d$ with log-characteristic function

 \beg{equation*}\beg{split} \log\hat\pi_t (u)&=  \int_0^t\eta(\e^{sA^*}u)\d s + tc_0S(|u|^2)\\
 &=\int_0^t \eta_2(\e^{sA^*}u)\d s
   +\int_0^t\big\{\eta_1(\e^{s A^*}u)\d s +c_0 S(|u^2|)\big\}\d s.\end{split}\end{equation*} Now, letting $P_t^S$ be the
   semigroup for the L\'evy process with L\'evy symbol $-c_0 S(|\cdot|^2)$, and letting

   $$ \tt P_t f(x)= \int_{\R^d}f(x+z)\pi_t(\d z),$$ we obtain from (\ref{M}), (\ref{M2}) and the definition of $\pi_t$ that

   $$ P_t f(x)= P_t^S \tt P_t f(\e^{tA}x).$$ Combining this with   (\ref{A1'}) we obtain

   $$\|\nn P_t f\|_\infty\le \|f\|_\infty \aa(c_0t).$$
   (V) For $t\in (0,1]$ and $r_0\in (0,\infty).$ Let $\rr_0, \bar L_t^0$ and $\bar L_t$ be in (II). Let

   \beg{equation*} \beg{split} &\bar P_t^1f(x)=\E\bigg\{f\bigg(\e^{tA} x +\int_0^t \e^{(t-s)A} \d\bar L_s^0\bigg)1_{\{\bar\tau_1\le t\}}\bigg\},\\
   &\bar P_t f(x)=\E f\bigg(\e^{tA} x +\int_0^t \e^{(t-s)A} \d\bar L_s\bigg).\end{split}\end{equation*}
Then (\ref{A4}) holds. Since (\ref{R}) holds for $\bar\nu$ in place of $\nu$, according to (IV) and the argument leading to (\ref{DDD}) using
Theorem \ref{T3.1},  there exists a constant $c\in (0,\infty)$ depending only on $d$ and $\theta$ such that

$$\|\nn \bar P_t\|_\infty \le\|f\|_\infty \aa(c_0t),\ \ \|\nn \bar P_t^1f\|_\infty \le \ff {c S(r_0^{-2}) t} {r_0}\|f\|_\infty.$$ Combining this
with (\ref{A4}) we derive the desired gradient estimate (\ref{G}).

\section{Proofs of Corollary \ref{C1.3} and Example \ref{E1.3}}

\beg{proof}[Proof of Corollary \ref{C1.3}]  Since the gradient estimate $\|\nn P_t^{+\si} f\|_\infty\le c(t)\|f\|_\infty$ is equivalent to

$$|P_t^{+\si} f(x)- P_t^{+\si} f(y)|\le c(t)\|f\|_\infty |x-y|,\ \ x,y\in\R^d,$$ by the monotone class theorem it suffices to prove for $f\in C_b^2(\R^d)$. By (\ref{GG}), in this case we have

$$\ff{\d}{\d s} P_s P_{t-s}^{+\si}f = P_s (\scr L-\scr L^{+\si}) P_{t-s}^{+\si} f = -P_s (\si P_{t-s}^{+\si} f),\ \ s\in [0,t].$$ Consequently,

$$P_t^{+\si} f= P_t f +\int_0^t P_s(\si P_{t-s}^{+\si} f) \d s.$$ Combining this with Theorem \ref{T1.1}, we finish the proof.\end{proof}

\beg{proof}[Proof of Example \ref{E1.3}]  (1) follows immediately from Theorem \ref{T1.1} and Corollary 1.2 by taking $S(r)= c r^{\aa/2}.$ To prove (2), we  take

$$S_\vv(r)= \log^{1+\vv}(1+r^{1/(1+\vv)}).$$ According to \cite{Aus}, for any Bernstein function $S$ and any $\dd>1$,
$r\mapsto S^\dd(r^{1/\dd})$ is again a Bernstein function. In this case we have

 $$\nu(\d z)\ge c 1_{\{|z|\le r_0\land 1\}} |z|^{-d} S_\vv(|z|^{-2}) \d z.$$  Then the desired gradient estimate follows immediately from Theorem \ref{T1.1}. \end{proof}

 \paragraph{Acknowledgement.} The author would like to thank Dr.  Jian Wang for a number of corrections and comments on   earlier versions of the paper, and to thank
 Professors M. R\"ockner, R. Schilling, R. Song, and   L. Wu for useful conversations. He also like to thank the referee for careful comments and    corrections. 
\beg{thebibliography}{99}

\bibitem{A} D. Applebaum, \emph{L\'evy Processes and Stochastic Calculus,} Cambridge University Press, 2004.

\bibitem{B} D. Bakry, \emph{On Sobolev and logarithmic Sobolev
inequalities for Markov semigroups,} ``New Trends in Stochastic
Analysis'' (Editors: K. D. Elworthy, S. Kusuoka, I. Shigekawa),
Singapore: World Scientific, 1997.

\bibitem{Bismut} J. M. Bismut, \emph{Large Deviations and the
Malliavin Calculus,} Boston: Birkh\"auser, MA, 1984.

\bibitem{BRS} V. I. Bogachev, M. R\"ockner, B. Schmuland,
\emph{Generalized Mehler semigroups and applications,} Probab. Th. Relat. Fields 105(1996), 193--225.

\bibitem{BJ} K. Bogdan, T. Jakubowski, \emph{Estimate of heat kernel of fractional Laplacian perturbed by gradient operators,} Comm. Math. Phys.
271(2007), 179--198.

\bibitem{EL}  K.D. Elworthy  and  X.-M. Li, \emph{Formulae for the
derivatives of heat semigroups,} J. Funct. Anal. 125(1994),
252--286.


\bibitem{GRW} M. Gordina, M. R\"ockner, F.-Y. Wang, \emph{Harnack
inequalities for subordinate semigroups,}  to appear in Potential Analysis. arXiv: 1004.3016.

\bibitem{KS} V. Knopova, R.L. Schilling, \emph{A note on the existence of transition probability densities for L\'evy processes,}
 arXiv:1003.1419.

 \bibitem{KS2} V. Knopova, R.L. Schilling, \emph{Transition  density estimates for a class of L\'evy and   L\'evy-type processes,}
 arXiv:0912.1482v1.

\bibitem{J} N. Jacob, \emph{Pseudo Differential Operators and Markov
Processes  (Volume I),}  Imperial College Press, London, 2001.

\bibitem{JS} T. Jakubowski, K. Szcypkowski, \emph{Time-dependent gradient perturbations of fractional Laplacian,}
J. Evol. Equ. 10(2010), 319--339.

\bibitem{PZ} E. Priola, J. Zabczyk, \emph{Densities for Ornstein-Uhlenbeck processes with jumps,} Bull. Lond. Math. Soc. 41(2009), 41--50.


\bibitem{ROW} M. R\"ockner, S.-X. Ouyang, F.-Y. Wang, \emph{Harnack inequalities and applications for Ornstein-Uhlenbeck semigroups with jump,} arXiv: 0908.2889.

\bibitem{RW03} M. R\"ockner, F.-Y. Wang, \emph{Harnack and functional inequalities for generalized Mehler semigroups,} J. Funct. Anal. 203(2003), 237--261.

\bibitem{Aus} L. R. Schilling, \emph{Subordination in the sense of Bochner and a related functional caculus,} J. Austral. Math. Soc. Ser. A 64(1998), 368--396.
\bibitem{SV} R. Song, Z. Vondra\v{c}ek, \emph{Potential theory of subordinate Brownian motion,} in $``$Potential Analysis of Stable Processes and its Extensions", 87-176, Lecture Notes in Math. Vol. 1980, Springer, Berlin, 2009.


\bibitem{W97} F.-Y. Wang, \emph{On estimation of the
logarithmic Sobolev constant and gradient estimates of heat
semigroups,} Probab. Theory Relat. Fields 108(1997), 87--101.

\bibitem{W10} F.-Y. Wang, \emph{Coupling for   Ornstein-Uhlenbeck   Processes with Jumps,} to appear in Bernoulli. \   arXiv: 1002.2890.

 \end{thebibliography}
\end{document}